# Locally controlled globally smooth ground surface reconstruction from terrestrial point clouds

Igor Rychkov

Department of Mathematics, University of Canterbury, New Zealand
`igor.rychkov@canterbury.ac.nz`

**Abstract.** Approaches to ground surface reconstruction from massive terrestrial point clouds are presented. Using a set of local least squares (LSQR) planes, the "holes" are filled either from the ground model of the next coarser level or by Hermite Radial Basis Functions (HRBF). Global curvature continuous as well as infinitely smooth ground surface models are obtained with Partition of Unity (PU) using either tensor product B-Splines or compactly supported exponential function. The resulting surface function has local control enabling fast evaluation.

**Keywords:** surface reconstruction, point clouds, least squares, Hermite radial basis functions, partition of unity, B-splines, curvature continuous

## 1 Introduction

In this paper we consider the problem of fitting a coarse level approximation to the geophysical data obtained with new surveying technology demanding better scalability of the modeling solutions. Terrestrial laser scanning (TLS) is a recent advancement in light detection and ranging (LiDAR) technology exploiting a new generation of high speed, high resolution laser scanners originally developed for structural engineering surveys. Applied to geomorphology [1] it can generate billions of points of terrestrial features in sub-centimeter detail and, by means of registering several overlapping point clouds, over kilometer long areas. The sheer size of these datasets demand a high performance processing and data organization to accomplish even a basic statistical analysis [1]. Information contained in the terrestrial point cloud allows multiscale study of morphological features down to the texture roughness of the surfaces. For example, by scanning fluvial sediment surfaces, one can map grain size distribution by calculating local standard deviations of the point displacements. The points in this case should be "de-trended" first by deducting a local ground profile. A ground model is therefore required. To design the criteria of an acceptable ground model is the first aim of this paper. Although we shall discuss the problem in the context of nearly horizontal, tilted river channel surfaces, the coarse modeling presented would be applicable to other terrestrial surfaces such as mountain slopes or even vertical cliffs after the "best" plane have been designated to be the $(x, y)$ grid plane.

In order to insure some consistency in the separation of the scales into "ground" and "sub-ground" scales the ground model should be at least continuous. Furthermore, there is a linked application field that may impose its own requirements on the surface model, that is the hydraulic modeling, for which the ground model should at least have no "holes" [2]. Finally, running a 2D CFD or particle fluid simulation on the curved ground surface of $z = f(x, y)$ (under development by the authors) would require it to have continuous second derivatives at least. Therefore, a curvature continuous ground model should be built in order to accommodate current and upcoming applications. The resulting ground model although globally smooth should still be locally controlled allowing a better scalability given the growing spans of the terrestrial surveys.

A key issue of terrestrial point clouds which is the point density inhomogeneity and particularly, the undersampled areas that we shall refer to as "holes". Holes arise from the physical limits of the laser beams that cannot "go" under water, for example, or simply from the lack of data points from a particular angle of some areas. Other areas can be declared holes if they are undersampled below certain threshold to ensure good statistical significance. One intuitive design principle regarding hole-filing is that the holes of certain scale should be interpolated over using the surrounding trends from the areas of *the same scale*. This rules out piecewise linear interpolation across empty areas using, e.g., Delaune triangulation.

Upon these considerations we suggest a three-stage approximation approach. Firstly, a set of local models is determined. Then the hierarchy of local models is built and uded to fill in the holes of corresponding scales. Alternatively, in Section 4 we introduce HRBF methodology at the hole-filling stage capable of inferring the local ground models by interpolation. In either case the final stage is to blend the local models together by a PU with tensor product cubic B-splines or $C^\infty$ compactly supported exponential as described in Section 3.

Applications of the approximation theory methods to the fitting of continuous global surfaces in geophysics have been addressed before, see for example [3] and references therein. The contribution of this paper is that while our surface is also globally smooth it is however locally controlled which could be an advantage for larger problems. Another novelty is an application of HRBF to benefit from the often available normal vector field and a couple of new PU blending functions.

## 2 Local ground model

A coarse, or "ground" model from is a much reduced representation of a high resolution terrestrial point cloud keeping the variations above certain "grid" scale. The higher frequency features are filtered out. This is the only stage where the whole, large point cloud is analyzed and the further modeling is based on the set of local ground models, a much smaller problem for typical values of the coarse model spacing and the range of scanning. For example, a point cloud over an $10^3 \times 10^2 m^2$ area with a $10^{-2}m$ resolution may contain a billion points, whereas the $1m$ grid will have only $10^5$ cells which is the maximum number of the points along with the attached slope normals defining the ground model at this scale.

We shall call our ground modeling locally controlled in contrast to what one might call a globally controlled ground modeling. The two approaches take the different answers to the question of whether the coarse model is to be determined from the local point positions or inferred from the surrounding trends disregarding the local points. The latter produces a ground surface model locally approximating the surface without finer scale local information. A multiresolution analysis using 2D tensor-product wavelets would be an example of such analysis, albeit inapplicable to the massive point clouds due to the high irregularity and the enormous size of the sampling data.

As the building blocks of our local ground model we chose local slopes defined as the best plane fitted by the (total) LSQR. A slope is therefore defined by the centroid and the normal vector. In order to define what the local means the $(x, y)$ plane is partitioned in a simple regular grid and the point coordinates are scaled so that the grid resolutions are always $(1, 1)$. In other words we will be working with centroid-thinned point clouds where the centroids also carry the normal vector.

It has to be noted that original points may not come from a functional surface $z = f(x, y)$. Essentially 3D features could be present, such as boulders with undercuts, overhanging cliffs and gorges, and vegetation. The LSQR planes simplify away all this complexity and, with a reasonable choice of the grid spacing, would not be too steep even when obtained by the total LSQR.

The locally controlled ground models use the original huge point cloud only once, a costly but straightforward operation [1]. Therefore, without loss of applicability to large datasets, in the algorithm development we chose to consider a minimal dataset consisting of merely 15666 points (Fig. 1) from the actual scan of a river bar of roughly 30 metres long. It has rather steep slopes and both holes to interpolate and outskirts to extrapolate. We refer to the both tasks as hole-filling below.

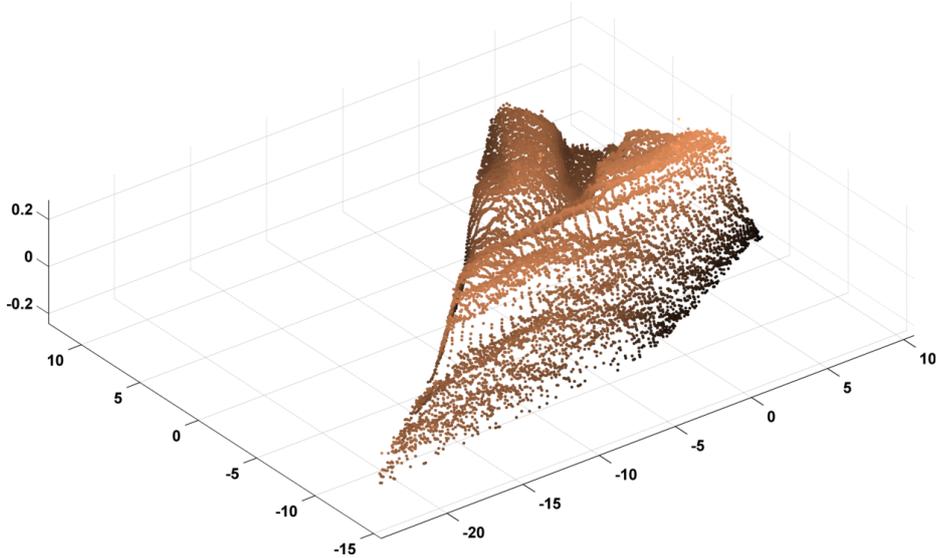

**Fig. 1.** The test point cloud dataset.

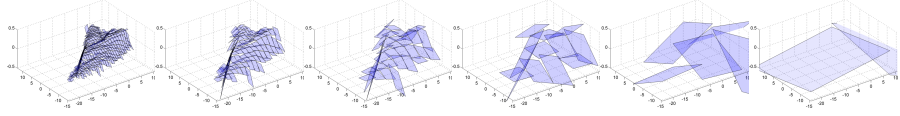

**Fig. 2.** Hierarchy of local slopes: slope grouping from left to right as the grid spacing doubles.

The first sub-plot in the Fig. 2 shows the slopes fitted in the data cells of the original grid as a set of rectangular patches, what may be called a "polygon soup". For the purposes of finding larger trends for hole-filling we combine the slopes into the next coarser grid, thus building a hierarchy of slopes until the first "full" grid with no holes is reached. Combining the slopes is done by re-fitting the LSQR plane to the slope vertices from each quadrant, maximum 16 points. The slopes at each stage appear to be a very crude approximation if they were to represent the local ground models. It prompts one to reconsider the questions: what is the ground profile and how to define it locally? The LSQR plane fitted through the data of a scanned tree, for example, would be very different from the soil profile on which it is planted. In fact, the local ground model should rather be defined as the slope with the tree removed, derived from the surrounding large-scale trends, analogous to the notion of the mean field in physics. Coming back to our situation, we define the local ground model as the 3x3 kernel mean filter of local slopes, i.e., slopes are replaced with their averages over the 3x3 neighborhood. Starting from the top full grid, ground models are projected onto the empty cells in the next finer sub-grid, averaging the slopes there, and so on to the original ground grid, Fig. 3.

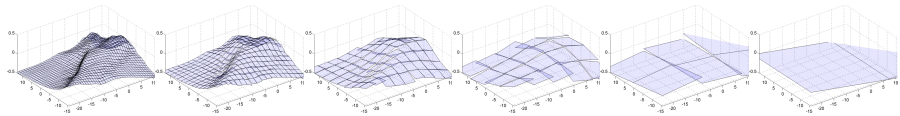

**Fig. 3.** Kernel smoothing and hole-filling: from right to left from the coarsest to the finest, original ground grid

## 3   PU blending splines

Slopes produced in this way look "almost" continuous and for some applications, e.g., roughness analysis by detrended standards deviation [?], such ground model would suffice and be even more convenient as the orthogonal projection of a point onto a plane is straightforward to find. One can go further and produce a globally $C^2$ smooth surface, also called curvature continuous surface, by blending the local slopes with a tensor product PU. We call a family of functions $\phi_i(x)$ a PU over given interval if $\sum_{i=1}^{\infty} \phi_i(x) = 1$ on the interval.

A suitable linear combination of cubic B-splines [4] should ensure both the $C^2$ smoothness and the PU over the central grid interval. Let $B_k^3$ be a cubic B-spline supported and positive on $x \in (t_k, t_{k+4})$, where $t_k$ are the nodes placed at the centers and the vertices of the grid intervals, so that $t_{2i}$ is the center of the grid interval $i$. For each interval $i$ we construct the following spline function, letting $k = 2i$,

$$\phi_i(x) = \frac{1}{2} B_{k-3}^3(x) + B_{k-2}^3(x) + \frac{1}{2} B_{k-1}^3(x) \tag{1}$$

Let us choose for the moment $t_0 = 0$ and look at $\phi_0(x)$ (Fig. 4). It is a symmetric 6-piecewise cubic polynomial overlapping three adjacent grid intervals, $[t_{2i-1}, t_{2i+1}], i = -1, 0, 1$. Therefore, the central interval is overlapped also by three functions only, $\phi_{i-1}, \phi_i, \phi_{i+1}$. Adding them together we recover the unity (Fig. 4) which is a consequence of the PU property of B-splines [4].

We have also constructed and tested a PU blending function that results in an infinitely differentiable $C^\infty$ global function,

$$\tilde{\phi}_i(x) = \begin{cases} 1 & x = 0 \\ 0 & |x| \geq a \\ \dfrac{1}{e^{s\left(\frac{1}{1-|\frac{x}{a}|} - \frac{1}{|\frac{x}{a}|}\right)} + 1} & \text{otherwise} \end{cases}, \tag{2}$$

which for the smoothing parameter $s = 1$, and the interval length $a = 2$ looks similar to Fig. 4 and also insures the PU over the central interval. However, it extends only to the centers of the left and right adjacent intervals (it is zero beyond the interval [-2,2]), and the overall blending is quite narrow. Paradoxically, this $C^\infty$ PU blender results in a less smooth surface overall when the local models are as crude as the LSQR planes. It would be appropriate for more precise local models that are in need of stitching smoothly at the narrow boundaries. Another argument in favor of the B-splines (1) is that the final approximate can be fully described locally by bi-variate polynomials further simplifying rapid, parallelizable evaluation.

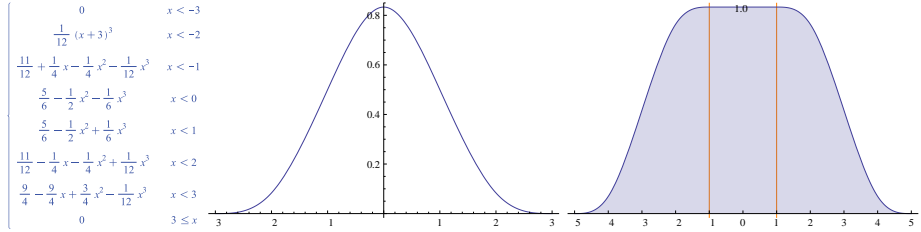

**Fig. 4.** Formulas (left) and a plot (center) of the spline function $\phi_0(x)$, and the unity sum $\phi_{-1} + \phi_0 + \phi_1$ (right)

For a grid cell with the duple index $(i, j)$ the tensor product bi-variate spline function is simply $\phi_i(x)\phi_j(y)$ having the corresponding tensor-product PU property. The tensor-product spline thus constructed can be used in blending the local ground approximants into a globally curvature continuous surface of a function,

$$g(x,y) = \sum_{\alpha=i-1}^{i+1} \sum_{\beta=j-1}^{j+1} \bar{s}_{\alpha\beta}(x,y)\phi_\alpha(x)\phi_\beta(y), \ \ (i,j) = (\lfloor x \rfloor, \lfloor y \rfloor) \quad (3)$$

Figure 5 shows the smooth ground function based on the averaged local slopes corresponding to the leftmost subplot in Fig. 3.

## 4 Filling in holes by HRBF interpolation

In Section 2 the holes in the ground model were filled in by slope grouping, averaging, and projecting, with the hope that the holes of a certain scale are filled in according to the trends from the surrounding areas of an appropriate scale. In the multivariate approximation theory the method of global radial basis functions is particularly known for its visually appealing and theoretically justified hole-filling properties [5]. Furthermore, since our slopes data contains not only heights but also normals, this can be related to the so called first order Hermite RBF problem [6] of fitting an RBF to the scattered measurements of both the function and its first derivatives.

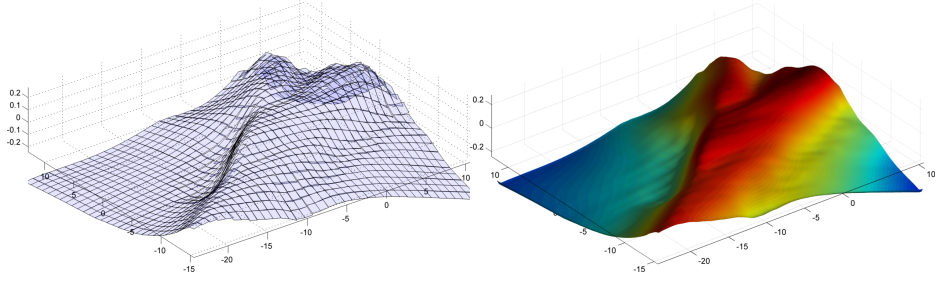

**Fig. 5.** Averaged local slopes (left) blended with the PU splines into a global $C^2$ surface (right)

We therefore look for the solution to the Hermit interpolation problem as follows. Given the values of the unknown bivariate function $f_j = f(\mathbf{x}_j)$ and its first derivatives $\nabla_j = (g_{2j-1}, g_{2j}) = \nabla f(\mathbf{x}_j) = (\partial_x f(\mathbf{x_j}), \partial_y f(\mathbf{x_j}))$ at $n$ scattered sites $\mathbf{x}_j$, find an interpolating function of the form

$$s(\mathbf{x}) = p(\mathbf{x}) + \sum_{j=1}^n c_j \psi(\mathbf{x} - \mathbf{x}_j) + \sum_{j=1}^n \mathbf{d}_j \cdot \nabla \psi(\mathbf{x}_j - \mathbf{x}), \qquad (4)$$

satisfying the interpolation conditions at the same locations $x_i$, $i = 1 \ldots n$,

$$s(\mathbf{x}_i) = f_i \qquad (5)$$
$$\nabla s(x_i) = \nabla_i, \qquad (6)$$

and the side condition for the bivariate polynomial $p$ of total degree $k - 1$,

$$\sum_l c_i p_l(\mathbf{x}_i) + \mathbf{d}_i \cdot \nabla p_l(\mathbf{x}_i) = 0, \qquad (7)$$

where $p_l$ is one of the $L = L(k-1)$ monomials of total degree not exceeding $k - 1$, and $k$ depends on the choice of the basic function $\psi$ [7]. In matrix form we have a linear system to solve,

$$\begin{bmatrix} A_{n \times n} & B_{n \times 2n}^T & Q_{n \times L} \\ B_{2n \times n} & C_{2n \times 2n} & R_{2n \times L} \\ Q_{L \times n}^T & R_{L \times 2n}^T & 0_{L \times L} \end{bmatrix} \begin{bmatrix} c_n \\ d_{2n} \\ a_L \end{bmatrix} = \begin{bmatrix} f_n \\ g_{2n} \\ 0_L \end{bmatrix}, \qquad (8)$$

with the elements $A_{ij} = \psi(\mathbf{x_i} - \mathbf{x}_j)$, $B_{ij} = \partial_i \psi(\mathbf{x}_{\lceil i/2 \rceil} - \mathbf{x}_j)$, $C_{ij} = -\partial_i \partial_j \psi(\mathbf{x}_{\lceil i/2 \rceil} - \mathbf{x}_{\lceil i/2 \rceil})$, $Q_{ij} = p_j(\mathbf{x}_i)$, $R_{ij} = \partial_i p_j(\mathbf{x}_i)$, where $\partial_i$ stands for $\partial_x$ when $i$ is odd and $\partial_y$ otherwise.

As the basic function $\psi$ we used Hardy's multiquadrics [8] $\psi(r) = \sqrt{r^2 + c^2}$ with $c = 0.1\delta$ where $\delta$ is an average spacing of the data locations, that is the scaled grid resolution in our case, $\delta = 1$. This value results in the lowest condition number of the matrix in (8), visually smooth interpolation and plausible extrapolation. A short review of the nature of "the notorious $c$ parameter" can be found on the p.237 of [9].

The method may not be suitable for larger grid sizes due to computational limitation on the size of the dense linear system (8). Fortunately, the number of equation that can be solved by modern out-of-core parallel solvers has reached 200 000 [10] which is of the order of a typical ground model grid for terrestrial scanning.

Solving the system (8) gives all the unknown coefficients in the RBF (4). The RBF can then be evaluated on the whole grid to fill in the holes with the interpolated slopes as shown in the left subplot of Fig. 6, where the rectangular patches representing the slopes are overlaid. In a straightforward implementation the evaluation of RBF functions may be rather costly an operation even for a centroid thinned point cloud. There are remedies to this difficulty. One is to implement a fast multipoles RBF method, as in [11,12]. Alternatively, at least for datasets which can fit in the GPU memory, a multi-point RBF evaluator can benefit from massively parallel execution model of recent GPU showing a 300-fold speedup as in our preliminary implementation using Nvidia's CUDA API.

Fortunately, for our purpose of hole-filling we only need to evaluate the RBF once to derive the local models on the whole grid. Then applying the PU spline blender from the Section 3 a surface quite similar to the original global HRBF is obtained (Fig. 6, right). The both approximations faithfully albeit slightly differently interpolate between individual slopes. The differences become even less significant when we apply a kernel smoothing to the underlying slopes in order to obtain the final curvature continuous ground model (7). Note that the overlaid plane patches are only nearly continuous, even though they may look as smooth as the surface mesh made of contour lines.

Although the surface obtained from the HRBF-derived slopes reproduces the essential features similarly to the surface obtained with hierarchical hole-filling shown in Fig. 5, the two clearly differ in the extrapolation "wings". This is not surprising given the *ad hoc* character of the "manual" hole filling in the latter, and some freedom in choosing the $c$ parameter of the multiquadrics RBF in the former. Ultimately, it should be left to the geomorphologist or hydrologist to decide on the extrapolation behavior of the ground surface model.

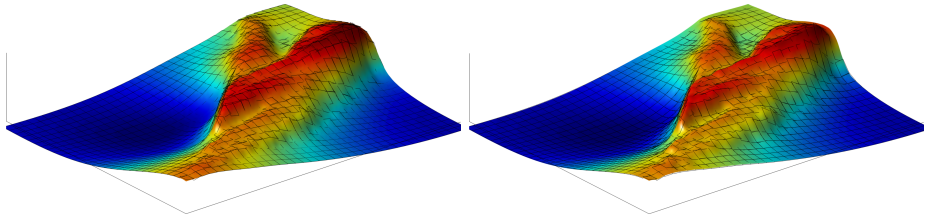

**Fig. 6.** HRBF fitted bumpy surface with re-evaluated slopes overlaid (left) and the same slopes blended PU into a $C^2$ bumpy surface

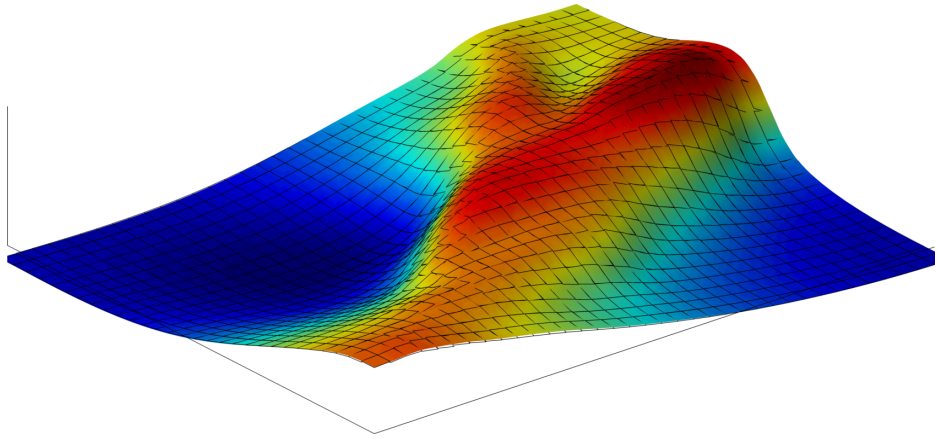

**Fig. 7.** Surface obtained from kernel averaged HRBF-derived slopes blended with the PU splines

# 5  Conclusion

We have presented approaches for ground model reconstruction from terrestrial point clouds. Sorting the point cloud into a regular 2d grid and fitting in each cell LSQR planes, we then examined two approached to fill in the holes, i.e., to impose such planes in the undersampled cells. In the first approach we continued to build the LSQR planes at the coarser grids until there were no holes. Using the 3x3 kernel averaging of the slopes to better define the smooth local ground approximation the grid was traversed hierarchy in the opposite direction filling in the holes from the information from the coarser level. We have also described an alternative approach solving a Hermite RBF interpolation problem with multiquadrics to fill in the holes of all scales "automatically". The two methods illustrate "local" and "global" approaches to hole-filling: combining slopes locally and interpolating from the overall set of slopes. Given the inherent arbitrariness in defining coarse models and in choosing the combining or fitting parameters in the two methods, the surfaces obtained are expected to agree only qualitatively and may have noticeably different extrapolation behavior. Having derived the local ground models everywhere on the grid, we used kernel-smoothing and blended them together into a global curvature continuous surface with a PU tensor product cubic B-splines. An infinitely smooth PU blender has also been presented constructed from smooth compactly supported exponential functions. Even without blending, the nearly-continuous set of the kernel-smoothed local ground slopes can be useful in some applications in fluvial geomorphology and sedimentology, while the $C^2$ ground surface can be used in hydraulics and 2D fluid modeling. The work allowed us to have a better understanding of coarse modeling in multi-stage and multi-scale surface reconstruction and compare RBF with *ad hoc* statistical fittings. It showed how to work with self-contradictory notion of locally defined and controlled coarse models despite its larger scale meaning. Work is underway by the authors to unlock the potential of our preferred method, the RBF, for multi-scale modeling, and for parallelizing RBF algorithms to make them more scalable.

**Acknowledgments.** We acknowledge a valuable discussion with Rick Beatson. We would also like to thank Damià Vericat and James Brasington for providing us with the TLS data and to acknowledge the interest and critique of participants in the Curves and Surfaces 2010 conference where this work was first presented.